\documentclass[12pt]{amsart}
\usepackage[mathscr]{eucal}
\usepackage{amsmath,amssymb,hyperref}
\usepackage{graphics}

\newtheorem{thm}{Th\'eor\`eme}[section]
\newtheorem{cor}[thm]{Corollaire}
\newtheorem{prop}[thm]{Proposition}
\newtheorem{lem}[thm]{Lemme}

\begin{document}

\title[Feuilletages en droites]
{Feuilletages en droites,\break \'equations des eikonales\break et autres \'equations diff\'erentielles}
\author[Dominique CERVEAU]{Dominique CERVEAU}
\address{\newline IRMAR, Campus de Beaulieu,\newline 35042 Rennes Cedex, France}
\email{dominique.cerveau@univ-rennes1.fr}
\date{}


\maketitle

\section*{Introduction}\label{S:intro}
Les solutions globales de certaines \'equations 
 aux d\'eriv\'ees partielles satisfont souvent des principes 
 de rigidit\'e \'etonnants. Ainsi un c\'el\`ebre th\'eor\`eme 
 de S. Bernstein affirme que si  $f : \mathbb R^2\longrightarrow \mathbb R$ 
 est de classe $C^2$ et a son graphe minimal, i.e $f$ satisfait l'\'equation 
 des surfaces minimales, alors $f$ est une fonction affine. Dans cet article 
 on s'int\'eresse en particulier aux solutions rationnelles 
 $f : \mathbb C^n\dashrightarrow \mathbb C$ 
 de l'\'equations des eikonales $E(f)= c^2$ o\`u $c\in\mathbb C$ et  
 $E$ d\'esigne l'op\'erateur
$$ E(f)=\sum \left(\frac{\partial f}{\partial z_{i}}\right)^2 . \leqno (1)$$
Cette \'equation est un cas sp\'ecial de l'\'equation diff\'erentielle : 
$$\hbox{\rm d\'et\ Hess}\ (f)=0 \leqno (2)$$
o\`u Hess $f$ d\'esigne la matrice Hessienne de $f$.

Si $f$ satisfait l'\'equation (2) l'application de Gauss :
$$G_{f} : z \rightsquigarrow \left(\frac{\partial f}{\partial z_{1}},\cdots, \frac{\partial f}{\partial z_{n}}\right) \leqno (3)$$
est d\'eg\'en\'er\'ee, i.e n'est pas de rang maximal. Si $f$ est rationnelle l'adh\'erence de l'image de $G_{f}$ est contenue dans une hypersurface alg\'ebrique. En particulier il existe un polyn\^ome $P$ irr\'eductible tel que
$$P\left(\frac{\partial f}{\partial z_{1}},\cdots, \frac{\partial f}{\partial z_{n}}\right)=0 \leqno (4)$$
et $f$ est solution d'une \'equation diff\'erentielle ne mettant en jeu que les d\'eriv\'ees partielles du premier ordre. Par exemple toute solution rationnelle $f$ de l'\'equation de Monge Amp\`ere~:
$$\frac{\partial^2 f}{\partial z_{1}^2}. \frac{\partial^2 f}{\partial z_{2}^2} - \left(\frac{\partial^2 f}{\partial z_{1} \partial z_{2}}\right)^2=1 \leqno (5)$$
satisfait une \'equation de type 4.

Dans le cas r\'eel il existe aussi un r\'esultat de rigidit\'e concernant les solutions globales de (5) ; plus pr\'ecis\'ement si 
$$f : \mathbb R^2 \longrightarrow \mathbb R$$
de classe $C^2$ est solution de l'\'equation de Monge-Amp\`ere (5), alors $f$ est un polyn\^ome de degr\'e plus petit ou \'egal \`a 2. Cet \'enonc\'e est attribu\'e \`a J¬\"orgen ; la d\'emonstration qui s'appuie sur le th\'eor\`eme de Picard (une fonction enti\`ere qui \'evite deux points est constante) permet de donner une preuve alternative \`a l'\'enonc\'e de Bernstein. Cette id\'ee est attribu\'ee \`a E. Bombieri.

On a longtemps cru avec Hesse qu'un polyn\^ome homog\`ene $f_{\nu}$ sur $\mathbb C^n$ satisfaisant
l'\'equation (2) d\'ependait en fait de moins de $n$ variables, cas o\`u l'on peut trouver un polyn\^ome $P$
satisfaisant (4) de degr\'e 1. Cet \'enonc\'e est correct en petite dimension $n\leq 4$ mais ne l'est plus
d\`es la dimension 5. L'exemple qui suit d\^u \`a Gordan et Noether \cite{GorNo} est tr\`es populaire car il intervient dans diff\'erents contextes \cite{DesCer}. Il s'agit de :
$$ \varphi (z) = z_{1}^2 z_{3}+ z_{1} z_{2} z_{4}+z_{2}^2 z_{5}. \leqno (6)$$
On remarque d'enbl\'ee que $\varphi$ est lin\'eaire dans les variables $z_{3}, z_{4} ,z_{5}$. Un calcul \'el\'ementaire montre que $\varphi$ satisfait l'\'equation diff\'erentielle : $$\frac{\partial\varphi}{\partial z_{3}} . \frac{\partial\varphi}{\partial z_{5}} -\left(\frac{\partial \varphi}{\partial z_{4}}\right)^2=0. \leqno (7)$$
Comme l'application de Gauss $G_{\varphi}$ est de rang maximal, l'adh\'erence de l'image de $G_{\varphi}$ est la quadrique de $\mathbb C^5$ d\'efinie par $X_{3} X_{5} -X_{4}^2=0$. Par suite il n'existe pas de polyn\^ome $P$ non constant de degr\'e 1 tel que l'on ait (4) : $P\left(\frac{\partial\varphi}{\partial z}\right)=0$.

Revenons \`a l'\'equation des eikonales ; elle trouve son importance dans les fondements de l'optique g\'eom\'etrique et aussi en m\'ecanique (Malus, Fresnel, Hamilton...). Si $f$ rationnelle satisfait $E(f)= c^2$ alors l'image de l'application de Gauss $G_{f}$ est contenue dans la quadrique $X^2_{1}+\cdots+X_{n}^2=c^2$. Consid\'erons, les coordonn\'ees $z_{1} 
,\cdots z_{n}$ \'etant fix\'ees (\`a l'action pr\`es du groupe engendr\'e par les translations et le groupe orthogonal complexe $O(n,\mathbb C)$) le champ de vecteur rationnel $X= \text{grad}\  f$ d\'efini par : 
$$X= \text{grad}\ f =\sum \frac{\partial f}{\partial z_{i}} \frac{\partial}{\partial z_{i}}. \leqno (8)$$
Comme on le sait au moins depuis Hamilton \cite{Ham}, si $\text{exp} t X$ est le flot (local, l\`a o\`u il a sens) du champ de vecteur $X$ alors
$$f\circ \text{exp}t X(x)= f(x)+ c^2 t \leqno (9)$$
 
\noindent (10)\quad et les trajectoires de $X$ sont contenues dans des droites (en fait sont d'adh\'erence des droites).

On doit imaginer, tout du moins en r\'eel, les niveau de $f$ comme un front d'onde (fig 1) se d\'epla\c{c}ant \`a la vitesse $c^2$ le long des droites param\'etr\'ees par :
$$t\rightsquigarrow z+t X(z). \leqno (11)$$

\begin{figure}[htbp]
\begin{center}

\input{artdom1.pstex_t}

\caption{}
\label{figure:1}
\end{center}
\end{figure}

On constate ainsi l'apparition naturelle de feuilletages alg\'ebriques de $\mathbb C^n$, et par suite de $\mathbb C \mathbb P(n)$, dont les feuilles sont d'adh\'erence des droites. Ce sujet \'etait tr\`es populaire chez les anciens g\'eom\`etres en particulier pour l'\'etude des surfaces et  l'int\'egration g\'eom\'etrique de certaines \'equations aux d\'eriv\'ees partielles. Dans ses "le\c{c}ons sur la th\'eorie g\'en\'erale des surfaces", Gaston Darboux y consacre un volume entier "les congruences de droites et les \'equations aux d\'eriv\'ees partielles". Il attribue \`a Malus le fait d'avoir le premier consid\'er\'e "de tels assemblages de droites".

En 1988 dans une courte note au CRAS \cite{Cer} j'ai propos\'e la classification des feuilletages (singuliers) en "droites" de
$\mathbb C^3$, r\'epondant en cel\`a \`a des pr\'eocupation de Ren\'e Thom et de l'astron\^ome Pecker. On trouvera en particulier dans cet article des r\'esultats annonc\'es dans cette note ainsi qu'une application \`a l'\'etude des solutions rationnelles de l'\'equations des eikonales.  Les \'enonc\'es pr\'ecis sont dans le chapitre qui suit.
\section{Champs de droites ou feuilletages en droites. Exemples et premiers r\'esultats.}
 Un feuilletage $\mathcal D$ en droites de l'espace $\mathbb C\mathbb P(n)$, ou de l'espace affine $\mathbb C^n$, est par d\'efinition un feuilletage alg\'ebrique (singulier) de dimension 1, tel qu'en tout point r\'egulier $m$ la feuille $\mathcal L_{m}$ passant par $m$ soit contenue dans une droite $D_{m}$. Si $\text{Sing} \mathcal D$ d\'esigne l'ensemble singulier de $\mathcal D$, $\text{Sing} \mathcal D$ est un ensemble de codimension sup\'erieure ou \'egale \`a deux ; visiblement $\mathcal L_{m}=D_{m}-\text{Sing} \mathcal D$ pour tout point r\'egulier $m$. Si $\mathbb C^n\subset\mathbb C\mathbb P(n)$ est une carte affine la restriction $\mathcal D_{/\mathbb C^n}$ est donn\'ee par un champ de vecteur polynomial :
 $$X=\displaystyle\mathop{\sum}_{i=1}^{n} X_{i}(z) \frac{\partial}{\partial z_{i}}$$
 o\`u les $X_{i}$ sont des polyn\^omes  tels que p.g.c.d $(X_{1},\cdots,X_{n})=1$. On a~ : 
 $$\text{Sing} \mathcal D \cap \mathbb C^n=\{z\in \mathbb C^n| X_{1}(z)=\cdots=X_{n}(z))=0\}.$$
 Si $m$ est un point r\'egulier de $\mathcal D$, par $m$ passe \'evidemment  une seule droite $D_{m}$ tangente \`a $\mathcal D$. A l'inverse, et nous le pr\'eciserons plus loin, si $m$ est un point singulier de $\mathcal D$ par $m$ passent une infinit\'e de droites qui sont, en dehors de $\text{Sing} \mathcal D$, des feuilles de $\mathcal D$. 
 
Le flot du champ $X$ satisfait l'\'equation diff\'erentielle 
$${\displaystyle\mathop{z}^{.}}_{i} (t) = X_{i} (z (t)),\quad z(0)=m. \leqno (12)$$
Dire que la trajectoire de $m$, param\'etr\'ee par $t\rightsquigarrow z(t)$, est rectiligne revient \`a dire que $\displaystyle\mathop{z}^{.} (t)$ et $\displaystyle\mathop{z}^{..} (t)$ sont colin\'eaires. Comme 
$$\begin{array} {ccc}
{\displaystyle\mathop{z}^{..}}_{i}(t)&=\displaystyle\mathop{\sum}^{n}_{k=1} \frac{\partial X_{i}}{\partial z k} (z (t)).{\displaystyle\mathop{z}^{.}}_{k} (t) &=\displaystyle\mathop{\sum}^{n}_{k=1} X_{k} (z(t)). \frac{\partial X_{k}}{\partial z k} (z (t))\\
 &=X(X_{i}) (z(t))\hfill&
\end{array} \leqno (13)$$ 
on constate que $X$ d\'efinit un feuilletage en droites si et seulement si :
$$\sum X(X_{i}) \frac{\partial}{\partial z_{i}} =  \mu . \displaystyle\mathop{\sum}^{n}_{i=1} X_{i} \frac{\partial}{\partial z_{i}} \leqno (14)$$
pour un certain polyn\^ome $\mu$, o\`u $X(X_{i})$ est la d\'eriv\'ee de $X_{i}$ le long de $X$.

En particulier les quotients $X_{i}/X_{j}$ sont des int\'egrales premi\`eres rationnelles du champ $X$.

A titre d'exemple consid\'erons dans $\mathbb C^2$ un feuilletage en droites d\'efini par le champ de vecteur local :
$$ X= X_{1} \frac{\partial}{\partial z_{1}}+ X_{2} \frac{\partial}{\partial z_{2}}.$$
Alors $a=X_{2}/X_{1}$ est solution de l'\'equation de Burger :
$$\frac{\partial a}{\partial z_{1}} + a. \frac{\partial a}{\partial z_{2}}= 0.$$
Cet exemple d'\'equation diff\'erentielle non lin\'eaire est bien connu des physiciens.

En dimension deux, dans $\mathbb C\mathbb P(2)$, les pinceaux de droites sont des exemples de champs de droites. Nous verrons plus loin que ce sont les seuls.

En dimension trois c'est un peu plus compliqu\'e. Dans $\mathbb C\mathbb P(3)$, consid\'erons un pinceau de plans donn\'e dans la carte affine $\mathbb C^3=\{(z_{1}, z_{2}, z_{3})\}$ par $z_{1}/z_{2}=t$. On se donne une application rationnelle du type :
$$t\rightsquigarrow (t z_{2}(t),  z_{2} (t), z_{3} (t))=m(t) \leqno (15)$$
Une  telle application revient \`a se donner dans chaque plan $z_{1}/z_{2}=t$ un point $m(t)$. Maintenant dans chacun de ces plans on consid\`ere le pinceau de droites de point de base $m(t)$. On feuillette ainsi $\mathbb C\mathbb P(3)$ (ou $\mathbb C^3$) en droites ; l'ensemble singulier est constitu\'e de l'union de l'axe des $z_{3}$ et de l'image de l'application (15), (fig 2) :

\begin{figure}[htbp]
\begin{center}

\input{artdom2.pstex_t}
 
\caption{}
\label{figure:2}
\end{center}
\end{figure}

Le premier cas est l'exemple g\'en\'erique.
Dans le second, l'image de l'application (15) coincide avec l'axe des $z_{3}$. Dans ces deux \'eventualit\'es nous dirons que nous avons un feuilletage radial dans les pages d'un livre ouvert. Dans le dernier l'application (15) est constante ; le champ de droites est "radial" au point $m=m(t)$.

Enfin introduisons la cubique gauche $\Sigma$ param\'etr\'ee dans une carte affine par 
$$t\rightsquigarrow (t, t^2, t^3)=\gamma (t). \leqno (16)$$
Pour chaque couple de points $\gamma (t_{1})$ et $\gamma (t_{2})$ de $\Sigma$ on m\`ene la s\'ecante \`a $\Sigma$ par ces deux points, convenant que si $t_{1}=t_{2}$ il s'agit de la tangente \`a $\Sigma$ en $\gamma (t_{1})$. On obtient ici encore un feuilletage de $\mathbb C \mathbb P(3)$ (ou $\mathbb C^3$) que nous appellerons feuilletage  associ\'e \`a la cubique gauche $\Sigma$ (fig 3).

\begin{figure}[htbp]
\begin{center}

\input{artdom3.pstex_t}
 
\caption{}
\label{figure:3}
\end{center}
\end{figure}

Parmi nos r\'esultats  en voici deux qu'il est ais\'e d'\'enoncer. Bien que leur preuves ne fassent intervenir que des arguments anciens et classiques, nous n'en avons pas trouv\'e trace.
\begin{thm}  Soit $\mathcal D$ un feuilletage en droites de $\mathbb C\mathbb P (3)$. Alors $\mathcal D$ est lin\'eairement conjugu\'e \`a l'un des exemples pr\'ec\'edents. Plus pr\'ecis\'ement $\mathcal D$ est de l'un des trois types suivants :

1. feuilletage radial en un point.

2. feuilletage radial dans les pages d'un livre ouvert.

3. feuilletage associ\'e aux cordes d'une cubique gauche.
\end{thm}
Nous \'etudierons quelques \'equations diff\'erentielle de type (2). En particulier en adaptant  le th\'eor\`eme   1.1 on obtiendra le :
\begin{thm}  Soit $f : \mathbb C^3 \dashrightarrow \mathbb C$ une solution rationnelle de l'\'equation aux eikonales $E(f)=c^2$, $c\in \mathbb C^*$. Il existe deux formes lin\'eaires $L_{1}= \alpha_{1} z_{1}+\alpha_{2} z_{2} +\alpha_{3} z_{3}$, $L_{2}=\beta_{1} z_{1}+\beta_{2}z_{2}+\beta_{3}z_{3}$ et $\ell \in \mathbb C(t)$ tels que :
$$f(z)=L_{1}+ \ell (L_{2})$$
avec $\|\alpha\|^2= \alpha_{1}^2+\alpha_{2}^2+\alpha_{3}^2=c^2$, $\|\beta\|^2=0$, $<\alpha |\beta>=\sum \alpha_{i} \beta_{i}=0$. Si $c=0$, alors $f$ est affine, $f=L_{1}$ avec $\|\alpha\|^2=0.$
\end{thm}
On en d\'eduit facilement le :
\begin{cor}  Soit $f : \mathbb R^3\dashrightarrow \mathbb R$ une solution rationnelle r\'eelle de l'\'equation aux eikonales. Alors $f$ est une fonction affine
\end{cor}
\section{G\'en\'eralit\'es sur les feuilletages en droites sur les espaces projectifs.}
Tous les champs de vecteurs holomorphes  rencontr\'es seront suppos\'es satisfaire $\text{cod}\ \text{Sing} X\geq 2$, condition \`a laquelle on se ram\`ene en divisant par le p.g.c.d des composantes de $X$.

Commen\c{c}ons par un \'enonc\'e facile en dimension 2 dont on retrouvera l'argument na¬\"if plusieurs fois.
\begin{prop}  Soit $\mathcal D$ un feuilletage en droites sur l'espace projectif $\mathbb C\mathbb P (2)$. Alors $\mathcal D$ est radial, ie correspond \`a un pinceau de droites concourantes. 
\end{prop}
\noindent\textbf{D\'emontration} : soit $\mathbb C^2 \subset \mathbb C\mathbb P(2)$. Si $\mathcal D_{/\mathbb C^2}$ n'a pas de singularit\'e alors les feuilles de $\mathcal D_{/\mathbb C^2}$ sont des droites parall\`eles. En particulier si $\mathcal D$ a une seule singularit\'e il s'agit d'un pinceau de droites. En g\'en\'eral si $D_{1}$ et $D_{2}$ sont deux droites tangentes \`a $\mathcal D$, elles se coupent en un point $M$. Un argument combinatoire imm\'ediat montre que toutes les autres droites passent par $M$.  \hfill$\square$

Nous donnons quelques \'enonc\'es pr\'ecisant le comportement des feuilles d'un champ de droites aux points singuliers.
\begin{lem}  Soit $X$ un champ de vecteurs holomorphes d\'efinissant un feuilletage (singulier) d'une  boule $B(0,\rho)$ dans $\mathbb C^n$. On suppose que chaque feuille r\'eguli\`ere $\mathcal L_{m}$ est contenue dans une droite $D_{m}$. Si $M$ est un point singulier de $X$, il existe au moins une droite $D$ passant par $M$ et tangente \`a $X$ ; \'eventuellement $D$ est contenue dans l'ensemble singulier de $X$.
\end{lem}
\noindent\textbf{Preuve} : on consid\`ere la vari\'et\'e d'incidence : 

$I:=\{(x,D), x\in \mathbb C\mathbb P (n)$, $D$ droite de $\mathbb C\mathbb P (n)\}$. Soit $m_{i}\in B(0,p)-\text{Sing} X$ une suite de points r\'eguliers convergant vers le point $M$. Comme la vari\'et\'e $I$ est compacte (on plonge $B(0,p) \subset \mathbb C^n$ dans $\mathbb C\mathbb P (n)$) la suite $(m_{i}, D_{mi})$ poss\`ede une sous suite convergente vers $(M,D)$. Par continuit\'e $X$ est tangent \`a $D$. \hfill$\square$

Nous pr\'ecisons le lemme pr\'ec\'edent :
\begin{lem}  Soit $X$ comme dans le lemme 2.2. On note $\mathcal C_{M}$ l'ensemble des droites $D$ passant par $M$ tel que $X$ soit tangent \`a $D$ (sur $B(0,p)$). Si le nombre de droites de $\mathcal C_{M}$ est fini alors $M$ est non singulier.
\end{lem}
\noindent\textbf{Preuve} : On choisit des coordonn\'ees  $z_{1}, \cdots z_{n}$ en $M$ telles que les droites $D_{1}\cdots,D_{s}$ de $\mathcal C_{M}$ ne soient pas contenues dans l'hyperplan horizontal $z_{n}=0$ et telles que $\text{Sing} X\cap (z_{n}=0)$ soit de codimension $\geq 3$. C'est possible. 

Ecrivons
$$X=\displaystyle\mathop{\sum}_{i=1}^n X_{i} \frac{\partial}{\partial z_{i}}.$$
Supposons que $X$ s'annule en $M=0$ ; alors l'ensemble $\sum=\{(X_{n}=0) \cap (z_{n}=0)\}$   est non trivial. On choisit une suite $m_{i} \in \sum$ et convergeant vers $M$. Alors les $D_{m_{i}}$ sont des droites horizontales et produisent des droites limites passant par $M$, horizontales et donc diff\'erentes des $D_{1},\cdots , D_{s}$. Ceci est absurde et par cons\'equent $X$ est non singulier en $M$. \hfill$\square$
\begin{cor} Si $M$ est un point singulier de $X$ l'ensemble $\mathcal C_{M}$ est un c\^one alg\'ebrique de dimension $\geq 2$.
\end{cor}
\noindent\textbf{Preuve} : $\mathcal C_{M}$ est un ensemble analytique ; comme c'est un c\^one, le th\'eor\`eme de Chow assure qu'il est alg\'ebrique. \hfill$\square$
\vskip 3mm
Voici encore une pr\'ecision qui sera utile pour d\'ecrire les singularit\'es de feuilletages en droites.
\begin{lem}  Soit $X$ holomorphe sur $B(0,\rho)$ et non singulier sur $B(0,\rho)-\{0\}$. On suppose que les feuilles r\'eguli\`eres $\mathcal L_{m}$ de $X$ sont contenues dans des droites $D_{m}$, $m\in B(0,\rho)-\{0\}$. Soit $D$ une droite tangente \`a $X$ passant par 0. On suppose que $X$ est non identiquement nul sur $D$. Si $D$ est isol\'ee dans le c\^one $\mathcal C_{0}$ (ie $D$ est une composante  irr\'eductible de $\mathcal C_{0}$) alors 0 est un point r\'egulier de $X$.
\end{lem}
\noindent\textbf{Preuve} : sans perdre de g\'en\'eralit\'e on suppose que $D$ est l'axe des $z_{n}$. On peut supposer aussi que $\rho>1$. Soit $H$ l'hyperplan $z_{n}=1$. Alors au voisinage du point $M_{0}=(0,\cdots,0,1)$, le champ $X$ est transverse \`a $H$. Par suite au voisinage de $M_{0}$ les feuilles sont param\'etr\'ees par les applications :
$$z'_{n} \rightsquigarrow (z_{1}+z'_{n} \eta_{1},\cdots, z_{n-1}+z'_{n}\eta_{n-1} , 1+ z'_{n})$$
o\`u $z_{n}=1+z'_{n}$, les $\eta_{i}$ sont holomorphes sur un voisinage de $M_{0}$ dans $H$ et satisfont $\eta_{i}(M_{0})=0$.

Evidemment les applications pr\'ec\'edentes sont globales en $z_{n}$. Comme la singularit\'e \'eventuelle de $X$ en 0 est isol\'ee dans $B(0,\rho)$ ainsi que dans $\mathcal C_{0}$, les droites ci dessus ne se coupent pas dans $B(0,\rho)$ et forment un voisinage de 0. D'autre part elles d\'efinissent visiblement un feuilletage r\'egulier au voisinage de 0. D'o\`u le lemme. \hfill$\square$

L'\'enonc\'e qui suit d\'ecrit les feuilletages en droites locaux \`a singularit\'e isol\'ee.
\begin{prop}  Soit $X$ un champ de vecteur holomorphe sur la boule $B(0,\rho)$ dans $\mathbb C^n$ \`a singularit\'e isol\'ee en 0. On suppose que les trajectoires de $X$ sont contenues dans des droites. Alors $X$ d\'efinit le feuilletage radial en 0, ie \`a unit\'e holomorphe multiplicative pr\`es $X=\displaystyle\mathop{\sum}_{i=1}^n z_{i} \frac{\partial}{\partial z_{i}}$.
\end{prop}
\noindent\textbf{D\'emonstration} : nous voulons montrer que chaque droite $D$ passant par 0 est tangente \`a $X$. Soit $D$ une telle droite que l'on suppose \^etre l'axe des $z_{1}$. Ecrivons :
$$X=\displaystyle\mathop{\sum}_{i=1}^n X_{i}\frac{\partial}{\partial z_{i}}$$
les $X_{i}$ \'etant holomorphes, $X_{i}(0)=0$.

Pla\c{c}ons nous sur l'ensemble analytique $\gamma$ :
$$\gamma:= \{X_{2}=\cdots=X_{n}=0\}.$$
Comme $X$ est \`a singularit\'e isol\'ee, $\gamma$ est une courbe passant par 0 sur laquelle $X_{1}$ ne s'annule qu'en 0 quitte \`a restreindre. Soit $m_{i}\in \gamma$ une suite tendant vers 0 ; visiblement avec les notations habituelles $D_{m_{i}}$ est parall\`ele \`a $D$ et n\'ecessairement $(m_{i}, D_{m_{i}})$ converge vers $(0,D)$. Par suite $X$ est tangent \`a $D$. \hfill$\square$
\vskip 3mm
\noindent\textbf{Remarques} : 1. on retrouve ainsi la preuve de la proposition 2.1.

2. Evidemment l'\'enonc\'e se globalise. Un feuilletage en droites de $\mathbb C^n$ ou $\mathbb C\mathbb P(n)$ ayant un point singulier isol\'e est radial.
\section{Classification des feuilletages en droites dans $\mathbb C\mathbb P (3)$}
Soit $\mathcal D$ un feuilletage en droite dans $\mathbb C\mathbb P(3)$. On peut supposer que $\mathcal D$ n'a pas de singularit\'e isol\'ee. Comme tout feuilletage de $\mathbb C\mathbb P(3)$ a des singularit\'es, l'ensemble singulier $\text{Sing} \mathcal D$ est de dimension pure 1 et est donc compos\'e de l'union $\Gamma_{1}\cup\cdots \cup \Gamma_{s}$ de courbes irr\'eductibles. Si $M\in \text{Sing} \mathcal D$,  on note encore $\mathcal C_{M}$ l'union des droites $D$ tangentes \`a $\mathcal D$ et qui passent par $M$. D'apr\`es le corollaire 2.4, $\mathcal C_{M}$ est un c\^one alg\'ebrique de dimension 2, avec \'eventuellement des branches de dimension 1. Mais le lemme 2.5 indique que l'ensemble :
$$\mathcal C^*_{M}=\overline{\mathcal C_{M}-\text{Sing} \mathcal D}$$
est une surface alg\'ebrique conique en $M$ ;
visiblement l'union :
$$\displaystyle\mathop{\cup}_{M\in\Gamma_{i}} \mathcal C_{M}^* \leqno (17)$$
est l'espace $\mathbb C\mathbb P(3)$ tout entier pour $i=1,\cdots,s$

Cette remarque implique que toute droite $D$ tangente \`a $\mathcal D$ coupe en au moins un point chaque composante $\Gamma_{i}$ du lieu singulier $\text{Sing} \mathcal D$. En particulier on obtient comme cons\'equence la 
\begin{prop}  Soit $\mathcal D$ un feuilletage en droites de $\mathbb C\mathbb P(3)$ dont l'ensemble singulier contient une droite $\mathcal D$. Alors $\mathcal D$ est radial dans les pages d'un livre ouvert.
\end{prop}
\noindent\textbf{D\'emonstration} : on consid\`ere le pinceau $P$ des plans  contenant $D$. Si $m\in \mathbb C\mathbb P(3)-\text{Sing} \mathcal D$ la droite $D_{m}$ coupe $D$. En particulier chaque plan $\pi$ de $P$ est $\mathcal D$ invariant et $\mathcal D_{/\pi}$ est un pinceau lin\'eaire de droites. \hfill$\square$
\vskip 3mm
\noindent\textbf{Remarque} : La description pr\'ecise de ce type de feuilletages, en particulier des singularit\'es, est donn\'ee en 1.
\vskip 3mm
Pour terminer la classification, on utilise avec les notations habituelles le :
\begin{lem} L'ensemble singulier $\text{Sing} \mathcal D$ a au plus 2 composantes.
\end{lem}
\noindent\textbf{Preuve} : soit $m$ un point r\'egulier et $D=D_{m}$ la droite tangente \`a $\mathcal D$ passant par $m$. Choisissons des coordonn\'ees affines $(z_{1}, z_{2}, z_{3})$ telles que $m=(0,0,0)$ et $D$ soit l'axe des $z_{3}$. Comme dans le lemme 2.5 nous param\'etrons les feuilles de $\mathcal D$ par les applications :
$$z_{3} \rightsquigarrow (z_{1}+z_{3} \eta_{1}, z_{2}+z_{3} \eta_{2}, z_{3})=F(z_{1},z_{2}, z_{3})$$
o\`u les $\eta_{i}$ sont des fonctions rationnelles en $(z_{1}, z_{2})$ r\'eguli\`eres en $(0,0)$.

Puisque $\mathcal D$ est r\'egulier en 0, l'application $F$ est un diff\'eomorphisme local en 0 ; en particulier le d\'eterminant Jacobien :
$$\hbox{\rm d\'e} JF=1+z_{3}\left(\frac{\partial \eta_{1}}{\partial z_{1}}+\frac{\partial \eta_{2}}{\partial z_{2}}\right) (z_{1}, z_{2})+z^2_{3}  \left(\frac{\partial \eta_{1}}{\partial z_{1}}.\frac{\partial \eta_{2}}{\partial z_{2}} -\frac{\partial \eta_{1}}{\partial z_{2}}. \frac{\partial \eta_{2}}{\partial z_{1}}\right) (z_{1}, z_{2})$$
est non nul en $(0,0,0)$. Ceci implique qu'\`a $(z_{1}, z_{2})$ fix\'es d\'et$JF$ ne s'annule qu'en deux valeurs de $z_{3}$ au plus. Comme toute droite de $\mathcal D$ coupe chaque composante du lieu singulier et qu'en chaque point singulier passent une infinit\'e de droites de $\mathcal D$, on en d\'eduit que $\text{Sing} \mathcal D$ a au plus deux composantes irr\'eductibles. \hfill$\square$
\vskip 3mm
Supposons que $\text{Sing} \mathcal D= \Gamma_{1}\cup \Gamma_{2}$. Soit $L\subset \mathbb C\mathbb P(3)$ une droite \'evitant $\text{Sing} \mathcal D$ et donc non tangente \`a $\mathcal D$. En particulier, pour tout point $m\in L$ la droite $D_{m}$ est transverse \`a $L$ et coupe $\Gamma_{i}$ en un seul point $M_{i}(m)$. Ceci d\'emontre que $\Gamma_{1}$ et $\Gamma_{2}$ sont rationnelles. 

D'apr\`es (17) l'ensemble des droites constituant $\mathcal D$ est pr\'ecis\'ement l'ensemble des droites joignant $\Gamma_{1}$ \`a $\Gamma_{2}$. En particulier les courbes $\Gamma_{1}$ et $\Gamma_{2}$ ne sont pas situ\'ees dans un m\^eme plan.
On remarque aussi que si $M_{1}$ est un point g\'en\'erique de $\Gamma_{1}$, le c\^one $\mathcal C_{M_{1}}$ contient $\Gamma_{2}$ : en effet, toute droite $D$ de $\mathcal C_{M_{1}}$ doit couper $\Gamma_{2}$.

Choisissons un 2-plan g\'en\'eral $\pi$ et soient 
$\{M_{1},\cdots,M_{t}\}=\Gamma_{1} \cap \pi$\quad , \quad $\{m_{1},\cdots,m_{s}\}=\Gamma_{2} \cap \pi$.

Les droites $[m_{i} M_{k}]\subset \pi$ joignant $m_{i}$ \`a $M_{k}$ sont dans $\mathcal D$. Mais sur chaque droite r\'eguli\`ere, il y a au plus 2 points singuliers. Par suite $s$ ou $t$ vaut 1, ie l'une des composantes $\Gamma_{i}$ est une droite ; cas d\'ecrit par la proposition 3.1. 

Dans la suite on suppose que $\text{Sing} \mathcal D$ se r\'eduit \`a une seule courbe irr\'eductible $\Gamma$ ; toujours d'apr\`es la proposition 3.1 on suppose encore que $\Gamma$ n'est pas une droite. Consid\'erons deux points distincts $M_{1}$ et $M_{2}$ de $\Gamma$.

Les c\^ones $\mathcal C_{M_{1}}$ et $\mathcal C_{M_{2}}$ se coupent le long d'une courbe qui est donc n\'ecessairement $\Gamma$. Par suite une droite g\'en\'erique de $\mathcal C_{M_{1}}$ coupe $\Gamma$ en $M_{1}$ et en un autre point de $\Gamma$. Les feuille de $\mathcal D$ sont donc exactement les cordes de $\Gamma$ (priv\'ees des singularit\'es) et les tangentes \`a $\Gamma$. Evidemment cel\`a implique que $\Gamma$ est une courbe gauche. Nous allons montrer que $\Gamma$ est une cubique gauche. Comme toujours choisissons $\pi$ un plan g\'en\'eral. Alors $\Gamma\cap \pi=\{M_{1},\cdots,M_{s}\}$ avec $s\geq 3$ puisque $\Gamma$ est irr\'eductible non plane. Les droites $D_{ij}$ joignant les $M_{i}$ \`a $M_{j}$ sont des droites de $\mathcal D$. Par connexit\'e de $\Gamma$, en chaque $M_{i}$ on a la m\^eme configuration pour les droites $D_{ij}$ ; d'autre part les droites $D_{ij}$ se coupent en des points de $\Gamma \cap M$. Comme sur chaque $D_{ij}$ on ne peut avoir que deux $M_{k}$, c'est \`a dire $M_{i}$ et $M_{j}$, n\'ecessairement $s=3$.

Maintenant \`a transformation lin\'eaire pr\`es il n'y a qu'une courbe gauche de degr\'e 3 dans $\mathbb C\mathbb P^3$ : la cubique gauche rationnelle $\sum$ param\'etr\'ee par (16). Nous avons ainsi d\'emontr\'e le th\'eor\`eme 1.1.
\section{L'\'equation d\'et Hess $f \equiv 0$.}
Consid\'erons une solution rationnelle $f : \mathbb C^n\dashrightarrow \mathbb C$  non triviale de l'\'equation (2) d\'et Hess $f \equiv 0$. Lorsque l'application de Gauss $G_{f}$ de $f$ est de rang g\'en\'erique $n-1$ nous dirons que $f$ est une solution maximale de (2). L'id\'eal $H(f)$ des polyn\^omes $Q\in \mathbb C[z_{1},\cdots,z_{n}]$ tels que $Q\left(\frac{\partial f}{\partial z}\right)$ est  alors engendr\'e par un polyn\^ome irr\'eductible $P$ ;  c'est le cas comme nous l'avons  vu pour le polyn\^ome de Gordan et N\"oether.

Les polyn\^omes $P$ que l'on peut ainsi obtenir ne sont certainement pas quelconques. Soient $f$ solution maximale de (2) et $P$ un polyn\^ome g\'en\'erateur de $H(f)$. La restriction de l'application de Gauss $G_{f}$ \`a un hyperplan g\'en\'eral param\`etre les z\'eros de $P$ : ainsi l'ensemble $(P=0)$ est une hypersurface unirationnelle.

Examinons plus pr\'ecis\'ement le cas de la dimension deux. Il existe alors une application rationnelle :
$$r :\mathbb C\dashrightarrow \mathbb C^2$$
$$t\rightsquigarrow r(t)=(r_{1}(t), r_{2}(t))\  ,\  r_{i}\in \mathbb C(t)$$
g\'en\'eriquement injective telle que $\overline{r(\mathbb C \mathbb P(1)}=\overline{(P=0)}\subset \mathbb C\mathbb P(2)$.

Par suite on dispose d'une factorisation :

$$\left\{
\begin{array}{cc}
\frac{\partial f}{\partial z_{1}} (z_{1},z_{2})=&r_{1} (\tau (z_{1},z_{2}))\\
\\
 \frac{\partial f}{\partial z_{2}} (z_{1},z_{2})=&r_{2} (\tau (z_{1},z_{2}))       
\end{array}
\right.\leqno (18)$$
o\`u $\tau : \mathbb C^2\dashrightarrow \mathbb C$ est rationnelle. 

Un calcul \'el\'ementaire montre que :
$$\frac{\partial r_{1}}{\partial t} (\tau (z)). \frac{\partial \tau}{\partial z_{2}}=\frac{\partial r_{2}}{\partial t}  (\tau (z)). \frac{\partial \tau}{\partial z_{1}}. \leqno (19)$$
En particulier le long d'un niveau $\tau=cste$,  la pente $\frac{\partial \tau}{\partial z_{2}}/ \frac{\partial \tau}{\partial z_{1}}$ est constante. Ainsi les niveaux de $\tau$ sont des droites, ce qui produit un feuilletage en droite de $\mathbb C^2\subset \mathbb C\mathbb P(2)$. Il y a deux cas, suivant que le point base de ce feuilletage soit ou non \`a distance finie, qui conduisent \`a $\tau$ de l'un des deux types :
$$\left\{
\begin{array}{cc}
\tau (z)=&\delta_{1} (a_{1}z_{1}+a_{2} z_{2}) \\
\tau (z)=&\delta_{2}\left(\frac{z_{1}-b_{1}}{z_{2}-b_{2}}\right)\hfill      
\end{array}
\right.\leqno (20)$$
o\`u les $a,b\in \mathbb C$ et $\delta_{i}\in \mathbb C (t)$.

Par une int\'egration \'el\'ementaire on obtient la : 
 \begin{prop}  Soit $f : \mathbb C^2 \dashrightarrow \mathbb C$ une solution rationnelle maximale de d\'et Hess $f=0$. Alors $f$ est de l'un des types suivants :
 
 1. $f=\ell_{1} (a_{1} z_{1} + a_{2} z_{2}) + c_{1} z_{1}+ c_{2} z_{2} + c_{3}$
 
 2. $f=\ell_{2} \left(\frac{z_{1}-b_{1}}{z_{2}-b_{2}}\right) . (z_{2}- b_{2})+ c_{1} z_{1} +c_{2} z_{2}+c_{3}$
 
 o\`u les $a_{i}, b_{i}, c_{i}$ sont des constantes et $\ell_{i} \in \mathbb C (t)$.
 \end{prop}
 On note que les solutions polynomiales de (2) sont de type 1. Dans le cas 1 le polyn\^ome $P$ g\'en\'erateur de l'id\'eal $H(f)$ est affine. Dans le cas 2 on peut tirer explicitement  en fonction de $f$ ce m\^eme polyn\^ome $P$. Il n'y a pas d'\'enonc\'e g\'en\'eral d\`es la dimension 3 d\'ecrivant les solutions de (2) ; toutefois dans la situation sp\'eciale o\`u l'on recherche les solutions de (2) sous forme polynomiales on dispose du :
\begin{thm}[\cite{BonEss}] Soit $f\in \mathbb C[z_{1}, z_{2}, z_{3}]$ un polyn\^ome solution de d\'et Hess $f=0$. Alors \`a conjugaison lin\'eaire pr\`es $f$ est de l'un des types suivants :

1. $\varphi_{1} = \varepsilon z_{1}+\varphi (z_{2}, z_{3}) \ ,\ \varphi \in \mathbb C [z_{2}, z_{3}]\ ,\ \varepsilon \in\{0,1\}$

2. $\varphi_{2}= a_{1}(z_{1})+z_{2} a_{2} (z_{1})+z_{3} a_{3} (z_{1}), a_{i} \in \mathbb C [z_{1}]$.
\end{thm}
Faisons quelques commentaires ; dans les deux cas  $f$ est affine dans une ou deux variables. Le graphe de l'application de Gauss de l'application $\varphi_{1}$ est contenu dans $z_{1}=\varepsilon$ ; si $\varphi_{1}$ est maximale le polyn\^ome $P$ est alors $z_{1}-\varepsilon$. Si $\varphi_{1}$ est non maximale il en est de m\^eme pour $\varphi$ ; dans ce cas $\varphi_{1}$ est du type $\varepsilon z_{1}+\varphi (z_{2})$ \`a conjugaison lin\'eaire pr\`es. Le graphe de l'application de Gauss est une droite. Dans le cas 2 en g\'en\'eral $\varphi_{2}$ est maximale. Si $P(z_{2}, z_{3})=0$ est l'\'equation de la courbe param\'etr\'ee par $t \rightsquigarrow (a_{2}(t), a_{3}(t))$ on a visiblement $P\left(\frac{\partial \varphi_{2}}{\partial z_{2}}, \frac{\partial \varphi_{2}}{\partial z_{3}}\right)=0$.
\section{L'\'equation des eikonales.}
Soit $f$ une fonction holomorphe d\'efinie sur un ouvert $V$ de $\mathbb C^n$ satisfaisant l'\'equation aux eikonales $E(f)=c^2, c\in \mathbb C$. Consid\'erons le champ de gradient de $f$ : 
$$X=\text{grad} f=\sum \frac{\partial f}{\partial z_{i}} \frac{\partial}{\partial z_{i}}, \leqno (8)$$
les coordonn\'ees $z_{i}$ \'etant fix\'ees.

Soit $t\rightsquigarrow z(t)$ une trajectoire de $X$ :
$${\displaystyle\mathop{z}^{.}}_{i} (t)=  \frac{\partial f}{\partial z_{i}} (z(t)). \leqno (21)$$ 
Suivant Hamilton consid\'erons l'acc\'el\'eration le long d'une trajectoire :
$${\displaystyle\mathop{z}^{..}}_{i} (t)=\displaystyle\mathop{\sum}_{k=1}^{n} \frac{\partial^2 f}{\partial z_{i} \partial z_{k}} (z(t)). {\displaystyle\mathop{z}^{.}}_{k} (t)= \displaystyle\mathop{\sum}_{k=1}^{n} \frac{\partial^2 f}{\partial z_{i}\partial z_{k}} (z(t)). \frac{\partial f}{\partial z_{k}} (z(t)). \leqno (22)$$
Si maintenant on d\'erive par rapport \`a $z_{i}$ l'\'equation $E(f)= c^2$, on obtient :
$$2  \displaystyle\mathop{\sum}_{k=1}^{n} \frac{\partial^2 f}{\partial z_{i}\partial z_{k}} . \frac{\partial f}{\partial z_{i}}=0. \leqno (23)$$
Ce qui implique la nullit\'e des ${\displaystyle\mathop{z}^{..}}_{i}$ ; en r\'esulte que, chaque fois qu'il est d\'efini, le flot $\varphi_{t}$ de $X$ est affine en $t$ ce qui indique en particulier que le feuilletage associ\'e \`a $X$ est un feuilletage en droites.

On remarque que, chaque fois que cel\`a a un sens : 
$$f\circ \varphi_{t} (z) = f(z) + c^2 t \leqno (24)$$
et :
$$\varphi_{t}(z)= z+ t. X(z)=z+t\frac{\partial f}{\partial z} (z) \leqno (25)$$
avec des notations \'evidentes. Ainsi les fonctions $\frac{\partial f}{\partial z_{i}} : V\rightarrow \mathbb C$ sont int\'egrales premi\`eres du champ $X$ : 
$$\frac{\partial f}{\partial z_{i}} (z+ t \frac{\partial f}{\partial z} (z))= \frac{\partial f}{\partial z_{i}} (z) \leqno (26)$$
autre traduction de l'\'egalit\'e (23).

Remarquons de suite que tout champ de droites n'est pas associ\'e \`a une solution de l'\'equation aux eikonales. Consid\'erons en effet le champ de droites dans $\mathbb C^3$ poss\'edant les int\'egrales premi\`eres $z_{3}$ et $\frac{z_{1}- r_{1} (z_{3})}{z_{2}-r_{2} (z_{3)}}$ o\`u $r_{1}$ et $r_{2}$ sont rationnelles. Alors si $X=\text{grad} f$ correspond \`a ce feuilletage on aura n\'ecessairement $\frac{\partial f}{\partial z_{3}}\equiv 0$ et   par suite $r_{1}$ et $r_{2}$ constants.

Nous allons maintenant, dans le cas de la dimension 3, utiliser le th\'eor\`eme 3 pour classifier les solutions rationnelles de l'\'equation aux eikonales.
\vskip 3mm
Soit $f=\frac{P}{Q}$, solution de $E(f)=c^2$ o\`u  $P$ et $Q$ sont des polyn\^omes sans facteur commun. 
\begin{lem}  Si $c\neq 0$, le feuilletage $\mathcal F$ en droites produit par $X= \text{grad} f$ n'est pas associ\'e aux cordes d'une cubique gauche.
\end{lem}
\noindent\textbf{D\'emonstration} : elle se fait par l'absurde. Remarquons que l'application de Gauss associ\'ee au feuilletage $\mathcal F$ par les cordes d'une cubique gauche est g\'en\'eriquement de rang $\geq 2$, et ceci pour tout choix de champ de vecteurs $Z$ rationnel d\'efinissant $\mathcal F$. En particulier l'image de l'application de Gauss :
$$G_{f} : z\sim \longrightarrow \left(\frac{\partial f}{\partial z_{1}} , \frac{\partial f}{\partial z_{2}}, \frac{\partial f}{\partial z_{3}}\right)$$
est dense dans la quadrique de $\mathbb C\mathbb P (3)$ donn\'ee par :
$$ z^2_{1}+ z^2_{2}+z^2_{3}=c^2 \leqno (27)$$
Soit $z_{0}$ un point g\'en\'erique de $\mathbb C^3$ o\`u $f=\frac{P}{Q}$ est holomorphe. On a l'\'egalit\'e entre fonctions rationnelles de $t$ :
$$\frac{P}{Q} \left(z_{0}+ t. \frac{\partial f}{\partial z} (z_{0})\right)=\frac{P}{Q} (z_{0})+c^2 t. \leqno (28)$$
Ecrivons $P$ et $Q$ sous forme d'une somme de polyn\^omes homog\`enes :
$$
\begin{array}{cc}
 P=& P_{0}+\cdots+P_{\nu}\\
 Q=&Q_{0}+\cdots+Q_{\mu}      
\end{array} \leqno (29)$$
les polyn\^omes $P_{\nu}$ et $Q_{\mu}$ \'etant non identiquement nuls. Remarquons, puisque $c\neq 0$, qu'un polyn\^ome homog\`ene non trivial ne peut s'annuler identiquement sur la quadrique (27) ; il existe donc un dense de $z_{0}$ pour lesquels on a $P_{\nu} \left(\frac{\partial f}{\partial z} (z_{0})\right)$ et $Q_{\mu} \left(\frac{\partial f}{\partial z} (z_{0})\right)$ non nuls. L'\'egalit\'e (28) se traduit alors au niveau des termes de plus haut degr\'e  par :
$$\cdots+t^{\nu} P_{\nu} \left(\frac{\partial f}{\partial z} (z_{0})\right)=\cdots +c^2 t^{\mu+1} Q_{\mu} \left(\frac{\partial f}{\partial z} (z_{0})\right)$$
qui implique $\nu=\mu+1$.

Puisque $G_{f}$ est dominante, le polyn\^ome $P_{\nu} - c^2 Q_{\nu-1}$ s'annule sur la quadrique (27). D'o\`u l'existence d'un polyn\^ome $K \in \mathbb C [z_{1}, z_{2}, z_{3}]$ tel que :
$$P_{\nu}-c^2 Q_{\nu-1}= (z^2_{1} +z^2_{2} +z^2_{3} - c^2). K \leqno (30)$$
Une fois encore on d\'eveloppe $K$ en somme de polyn\^omes homog\`enes :
$$K= K_{\alpha}+\cdots+K_{\beta}\quad,\quad \alpha\leq \beta.$$
et l'on observe en calculant les termes de plus haut et plus bas degr\'e de (30) que $\alpha=\nu-1$ et $\nu=\beta+2$. Ce qui est absurde. \hfill$\square$

Dans le lemme qui suit on traite le cas o\`u $c=0$ en utilisant une approche plus g\'eom\'etrique.
\begin{lem}  Le feuilletage en droites associ\'e \`a une solution rationnelle de l'\'equation aux eikonales $E(f)=0$ n'est pas du type cordes d'une cubique gauche.
\end{lem}
\noindent\textbf{D\'emonstration} : sous les hypoth\`eses du lemme $f$ est int\'egrale premi\`ere de son gradient $X=\text{grad} f$. Supposons que le feuilletage associ\'e \`a $X$ ait pour trajectoires g\'en\'eriques les cordes d'une cubique gauche $\Gamma$. Comme on l'a vu les $\frac{\partial f}{\partial z_{i}}$ sont aussi int\'egrales premi\`eres du champ $X$. Par suite les fibres  $f^{-1}(c)$ de $f$ sont des surfaces r\'egl\'ees par les trajectoires de $X$ et le long de ces trajectoires le plan tangent \`a $f^{-1}(c)$ est "constant". En r\'esulte que les $f^{-1} (c)$ sont des c\^one invariants par $X$, et donc les c\^ones $\mathcal C_{M}, M\in \Gamma$. Ceci vient du fait que les cordes de la cubique $\Gamma$ feuillettent $\mathbb C\mathbb P (3)-\Gamma$ et ne peuvent donc se rencontrer en dehors de $\Gamma$. Mais la famille des c\^ones $\mathcal C_{M}$ ne feuillettent pas $\mathbb C\mathbb P (3)-\Gamma$ ; il s'agit en fait d'un bi-feuilletage et par chaque point de $\mathbb C\mathbb P(3)-\Gamma$ passent deux tels c\^ones. D'o\`u une contradiction. \hfill$\square$
\vskip 3mm
Supposons que le champ $X= \text{grad} f$ soit tangent \`a un pinceau de plans $L_{1}/L_{2}=cst$ o\`u les $L_{i}$ sont affines non constants. Quitte \`a faire une translation on peut supposer les $L_{i}$ lin\'eaires et l'on \'ecrira :
$$
\begin{array}{cc}
 L_{1}=& \alpha_{1} z_{1}+\alpha_{2}z_{2}+\alpha_{3}z_{3}\\
 L_{2}=&\beta_{1} z_{1}+\beta_{2}z_{2}+\beta_{3}z_{3}. 
\end{array} \leqno (31)$$
On pose $<\alpha |\beta>= \sum \alpha_{i} \beta_{i}, \|\alpha\|^2=\sum \alpha^2_{i}$ et $\|\beta\|^2=\sum \beta^2_{i}$, En \'ecrivant explicitement que $L_{1}/L_{2}$ est int\'egrale premi\`ere du champ $X$, on obtient~ :
$$\left\{(\beta_{1}L_{1}-\alpha_{1}L_{2})\frac{\partial}{\partial z_{1}}+(\beta_{2}L_{1}-\alpha_{2} L_{2}) \frac{\partial}{\partial z_{2}}+(\beta_{3}L_{1}-\alpha_{3}L_{2}) \frac{\partial}{\partial z_{3}}\right\} (f)=Y(f)=0. \leqno (32)$$
Ainsi le champ lin\'eaire $Y$ s'annule sur la droite $L_{1}=L_{2}=0$ et poss\`ede $f$ pour int\'egrale premi\`ere. Le fait que ses composantes soient li\'ees produit une forme lin\'eaire non triviale :
$$L_{3}= \gamma_{1}z_{1}+\gamma_{2}z_{2}+\gamma_{3} z_{3} \leqno (33)$$
telle que $Y(L_{3})=0$. Remarquant que 
$$Y(L_{3})= <\beta | \gamma> L_{1}+<\alpha | \gamma> L_{2}$$
on en d\'eduit que :
$$<\alpha | \gamma>=<\beta | \gamma>=0.$$
Remarquons aussi que :
$$\left\{
\begin{array}{cc}
 Y(L_{1})=& <\alpha | \beta>L_{1} -\|\alpha\|^2 L_{2}\\
 Y (L_{2})=&\|\beta\|^2 L_{1} -<\alpha | \beta> L_{2}.
\end{array}\right. \leqno (34)$$
Ainsi $Y$ agit lin\'eairement sur l'espace vectoriel Vect$_{\mathbb C} (L_{1},L_{2})$ ; comme cette action est de trace nulle  on peut choisir $L_{1}$ et $L_{2}$ de sorte que :
$$\left\{
\begin{array}{cc}
 Y(L_{1})=& \lambda L_{2}\\
 Y (L_{2})=&-\lambda L_{2}.
\end{array}\right. \leqno (35)$$
o\`u $\lambda \in \mathbb C$. Ceci revient donc \`a supposer que $\|\alpha\|^2=\|\beta\|^2=0$ et $<\alpha | \beta >=\lambda$. On note alors, puisque $L_{1}$ et $L_{2}$ sont ind\'ependants, que $\lambda$ est non nul. En particulier le corps des int\'egrales premi\`eres rationnelles de $Y$ est engendr\'e par $L_{3}$ et la forme quadratique $L_{1} L_{2}$ ;  on note aussi que $L_{1}, L_{2}$ et $L_{3}$ sont ind\'ependants $(\lambda \neq 0)$. Par cons\'equent les trajectoires de $Y$ sont les coniques :
$$\left\{
\begin{array}{cc}
 \hfill L_{3}&= \text{cste}\\
L_{1} L_{2}&=\text{cste}.
\end{array}\right. \leqno (36)$$
Le feuilletage $\mathcal F$ associ\'e \`a $f$ (ses feuilles sont les fibres de $f$) est invariant par $Y$ ; s'il est d\'efini par la 1-forme polynomiale
$$\omega=A_{1}dz_{1}+A_{2}dz_{2}+A_{3} dz_{3} \leqno (37)$$
son lieu singulier :
$$\text{Sing} \mathcal F=\{A_{1}=A_{2}=A_{3}=0\} \leqno (38)$$
est aussi invariant par $Y$, et donc form\'e de trajectoires de $Y$. Le feuilletage associ\'e \`a $X=\text{grad} f$ est donn\'e par le champ $Z$
$$Z= A_{1} \frac{\partial}{\partial z_{1}} + A_{2} \frac{\partial}{\partial z_{2}} + A_{3} \frac{\partial}{\partial z_{3}}. \leqno (39)$$
Mais d'apr\`es le th\'eor\`eme 1.1 son lieu singulier est constitu\'e de la droite $L_{1}=L_{2}=0$ et \'eventuellement d'une courbe rationnelle $\Gamma$ coupant le plan g\'en\'erique $L_{1}= t L_{2}$ en un seul point en dehors de l'axe $L_{1}=L_{2}=0$.
Visiblement $\Gamma$ ne peut-\^etre une fibre g\'en\'erique de $(L_{1} L_{2}, L_{3})$ qui coupe deux fois chaque $L_{1}= t L_{2}$ en dehors de l'axe $L_{1}=L_{2}=0$ ; ni une fibre sp\'ecial (contenue dans $L_{1} L_{2}=0)$ qui n'appara\^{\i}t pas comme lieu singulier de feuilletages en droites. Ne reste que la possibilit\'e o\`u le lieu  singulier du champ de droites $Z$ est pr\'ecis\'ement l'axe $L_{1}= L_{2}=0$

\begin{figure}[htbp]
\begin{center}

\input{artdom4.pstex_t}
 
\caption{}
\label{figure:4}
\end{center}
\end{figure}

On sait que dans cette situation le champ $X=\text{grad} f$ poss\`ede les deux int\'egrales premi\`eres $L_{1}/L_{2}$ et $\frac{L_{3}- r(L_{1}/L_{2})}{L_{1}}$ o\`u $r\in \mathbb C (t)$ est une certaine fonction rationnelle. Comme $f$ est int\'egrale premi\`ere du champ $Y$, il existe une fonction rationnelle $\varphi \in \mathbb C (t_{1}, t_{2})$ telle que :
$$f=\varphi (L_{1}L_{2}, L_{3}). \leqno (41)$$
De sorte que $X = \text{grad} f$ s'\'ecrit :
$$X=\frac{\partial \varphi}{\partial t_{1}}\{ L_{1} \text{grad} L_{2}+ L_{2} \text{grad} L_{1}\}+\frac{\partial \varphi}{\partial t_{2}}  \text{grad} L_{3}. \leqno (42) $$
En \'ecrivant explicitement que $\frac{L_{3}-r(L_{1}/L_{2})}{L_{1}}$ est int\'egrale premi\`ere du champ $X$ on obtient :
$$<\alpha |\beta> . \frac{\partial \varphi}{\partial u} (L_{1} L_{2}, L_{3}). \left(L_{3}-r\left(\frac{L_{1}}{L_{2}}\right)\right)-\|\gamma\|^2 \frac{\partial \varphi}{\partial v}(L_{1} L_{2}, L_{3})=0.$$
Comme $<\alpha | \beta>=\gamma$ est non nul et comme $L_{1}, L_{2}, L_{3}$ sont ind\'ependants n\'ecessairement $\|\gamma\|^2\neq 0$ (utiliser $<\gamma | \beta>=<\gamma | \alpha >=0)$. Finalement $L_{3}-r (\frac{L_{1}}{L_{2}})$ appara\^{\i}t comme fonction de $L_{3}$ et de la forme quadratique $L_{1} L_{2}$ ; ce qui est absurde sauf si $r$ est constante. Mais dans ce cas le champ $\text{grad} f$ est radial.

Lorsque $\text{grad} f$ est radial :
$$\text{grad} f=h. \displaystyle\mathop{\sum}_{i=1}^3 z_{i} \frac{\partial}{\partial z_{i}} \quad , h\in \mathbb C (z_{1}, z_{2}, z_{3})$$
on a : 
$$ z_{i} \frac{\partial h}{\partial z_{j}}- z_{j} \frac{\partial h}{\partial z_{i}}=0. \leqno (43)$$
En particulier il existe $\ell \in \mathbb C(t)$ telle que :
$$h=\ell (z^2_{1}+ z^2_{2} +z^2_{3})$$
Par suite $f=L (z^2_{1}+ z^2_{2} +z^2_{3})$. L'\'equation des eikonales se traduit par 
$$4 t L'(t)^2=c^2$$
qui s'int\`egre en $L(t)=\pm c\sqrt t + cste$. La solution produite n'est pas rationnelle comme on le voit et comme le savait Malus... 

Il reste finalement un cas \`a examiner, celui o\`u le champ de droite $X=\text{grad} f$ est tangent \`a un pinceaux d'hyperplans parall\`eles $L_{3}=\text{cste}$ ; on pose 
$$L_{3}=\gamma_{1} z_{1} +\gamma_{2} z_{2} +\gamma_{3} z_{3}.$$
Les formes lin\'eaires $a_{1} z_{1} +a_{2} z_{2} + a_{3} z_{3}$ avec $\|a\|^2=c^2$ sont solutions de $E(f)=c^2$ et font partie du cas pr\'ec\'edent.

Puisque :
$$\gamma_{1} \frac{\partial f}{\partial z_{1}} +\gamma_{2} \frac{\partial f}{\partial z_{2}} + \gamma_{3} \frac{\partial f}{\partial z_{3}} =0$$
le champ de vecteur $Y=\sum \gamma_{i} \frac{\partial}{\partial z_{i}}$ annule $f$. Consid\'erons deux formes lin\'eaires ind\'ependantes 
\begin{eqnarray*}
L_{1}&=&\alpha_{1} z_{1} +\alpha_{2} z_{2} +\alpha_{3} z_{3}\\
L_{2}&=&\beta_{1} z_{1} +\beta_{2} z_{2} + \beta_{3} z_{3}
\end{eqnarray*}
telles que $<\alpha | \gamma>=<\beta | \gamma>=0$. Ces deux formes lin\'eaires engendrent le corps des int\'egrales premi\`eres de $Y$ si bien que 
$$f=\varphi (L_{1}, L_{2}) \quad \hbox{\rm o\`u}\quad \varphi \in \mathbb C (t_{1}, t_{2}). \leqno (44)$$
L'\'equation aux eikonales $E(f)=c^2$ se traduit par :
$$\|\alpha\|^2 \left(\frac{\partial \varphi}{\partial t_{1}}\right)^2+ \|\beta\|^2 \left(\frac{\partial \varphi}{\partial t_{2}}\right)^2 +2<\alpha | \beta> \frac{\partial \varphi}{\partial t_{1}}\  \frac{\partial \varphi}{\partial t_{2}}=c^2.\leqno(45)$$
Remarquons que $\|\alpha\|^2, \|\beta\|^2$ et $<\alpha, \beta>$ ne peuvent \^etre simultan\'ement nuls puisque $L_{1}$ et $L_{2}$ sont ind\'ependantes. Par suite $\varphi \in \mathbb C (t_{1}, t_{2})$ satisfait : 
$$\hbox{\rm d\'et\ Hess} \varphi\equiv 0.$$
La proposition 4.1 assure qu'il existe deux formes lin\'eaires $u_{1}$ et $u_{2}$ en deux variables telles que $\varphi$ soit de l'un des deux types :

$\varphi_{1}=\varepsilon  u_{1}+ \ell (u_{2})$

$\varphi_{2}=\varepsilon u_{1}+ \ell \left(\frac{u_{1}}{u_{2}}\right) u_{2} \quad ,\quad \ell \in \mathbb C (t), \varepsilon \in \{0,1\}.$

\noindent Quitte \`a effectuer un changement de notation on peut donc supposer que $f$ est de l'un des deux types :
$$f_{1}= \varepsilon L_{1} + \ell (L_{2}) \leqno (46)$$
$$f_{2}=\varepsilon L_{1}+\ell \left(\frac{L_{1}}{L_{2}}\right). L_{2}\qquad \varepsilon \in \{0,1\}. \leqno(47)$$
Dans le cas (46) l'\'equation des eikonales se traduit par :
$$\varepsilon \|\alpha\|^2+\ell' (L_{2})^2 \|\beta\|^2 +2\ell' (L_{2})<\alpha | \beta>=c^2. \leqno (48)$$
En particulier si $\ell'$ est non constante on aura $<\alpha |\beta>=\|\beta\|^2=0$ et $\varepsilon\|\alpha\|^2=c^2$ ; on peut alors, si $c\neq 0$, supposer que $\varepsilon =1$. Si $c=0$ alors $f$ est affine. Dans le cas (47) l'\'equation des eikonales conduit \`a :
$$\|\alpha\|^2 (\varepsilon +\ell' (t))^2+\|\beta\|^2 (\ell (t)- t \ell' (t))^2+2 <\alpha | \beta> (\varepsilon+ \ell'(t)) (\ell (t)-t \ell' (t))=c^2. \leqno (48)$$
Elle ne produit pas de nouvelles solutions rationnelles comme nous allons le voir en l'int\'egrant explicitement. En d\'erivant (48) on obtient :
$$\ell'' (t)\{ \ell' (t) (\|\alpha\|^2-2<\alpha | \beta> t-t^2 \|\beta\|^2)+\leqno (49)$$
$$\ell (t) (t\|\beta\|^2+<\alpha | \beta >)+\|\alpha\|^2 \varepsilon - \varepsilon t <\alpha | \beta>\}=0.$$
Evidemment (49) poss\`ede des solutions affines mais qui sont de type (46). Apr\`es simplification par $\ell''$ on obtient une \'equation diff\'erentielle lin\'eaire avec second membre. On remarque que $t\rightsquigarrow -\varepsilon t$ en est une solution particuli\`ere. L'\'equation sans second membre :
$$y' (\|\alpha\|^2- 2<\alpha | \beta>t- t^2 \|\beta\|^2)+y (t\|\beta\|^2+<\alpha | \beta>=0 \leqno (50)$$
s'int\`egre  en :
$$y(t) = cste (t^2 \|\beta\|^2+2<\alpha | \beta> t-\|\alpha\|^2)^{\frac{1}{2}} \leqno (51)$$
qui produit donc des solutions explicites de (48). Ces solutions ne sont pas rationnelles sauf lorsque $t^2 \|\beta\|^2+2<\alpha |\beta> t-\|\alpha\|^2$ est un carr\'e ~; ce cas conduit encore \`a une solution de type 46.

On obtient in fine le :
\begin{thm}  Les solutions rationnelles de $E(f)=c^2$, $c\neq 0$ sont de type $L_{1}+\ell (L_{2})$ o\`u les $L_{i}$  sont des formes lin\'eaires, $L_{1}=\alpha_{1} z_{1}+ \alpha_{2} z_{2} +\alpha_{3} z_{3}, L_{2}=\beta_{1}z_{1} +\beta_{2}z_{2}+\beta_{3}z_{3}$, satisfaisant $\|\alpha\|^2=c^2$, $\|\beta\|^2=0$, $<\alpha |\beta>=0$ ; $\ell$ est une fonction rationnelle. Les solutions rationnelles de $E(f)=0$ sont affines.
\end{thm}
\section{Automorphismes des feuilletages en droites.}
Soit $\mathcal F$ un feuilletage en droites de $\mathbb C\mathbb P (3)$ ; nous allons d\'ecrire quelques  groupes $Aut \mathcal F$ o\`u 
$$\text{Aut} \mathcal F :=\{\varphi \in \text{Aut} \mathbb C \mathbb P (3)\ \  | \ \ \varphi^* \mathcal F=\mathcal F\}.$$
Rappelons que $\text{Aut} \mathbb C \mathbb P (3) \cong P G L (4,\mathbb C)$ ; visiblement $\text{Aut} \mathcal F$ est un sous groupe alg\'ebrique de $\text{Aut} \mathcal F$. Si $\varphi$ est un \'el\'ement de $\text{Aut} \mathcal F$ alors $\varphi (\text{Sing} \mathcal F)=\text{Sing} \mathcal F$ ; dit autrement $\text{Aut} \mathcal F$ est un sous groupe de $\text{Aut}(\text{Sing} \mathcal F)$ le groupe des automorphismes de $\mathbb C\mathbb P(3)$ qui pr\'eservent $\text{Sing} \mathcal F$.
\begin{prop} soit $\mathcal F$ un feuilletage en droites associ\'e \`a une cubique gauche $\Gamma$ ; alors :
$$\text{Aut} \mathcal F =\text{Aut} (\text{Sing} \mathcal F)=\text{Aut} \Gamma \cong P G L(2,\mathbb C).$$
\end{prop}
\noindent\textbf{D\'emonstration} : si $\sigma : \mathbb C\mathbb P(1) \longrightarrow \mathbb C\mathbb P (3)$ est une param\'etrisation de $\Gamma$, alors chaque \'el\'ement de $\text{Aut} \mathbb C\mathbb P(1)\cong P G L (2,\mathbb C)$ se rel\`eve \`a $\Gamma$ et s'\'etend en un automorphisme de $\mathbb C\mathbb P(3)$. Comme un \'el\'ement $\varphi \in \text{Aut} \mathbb C\mathbb P(3)$ dont la restriction \`a $\Gamma$ est l'identit\'e en lui m\^eme, on a  $\text{Aut} \Gamma \cong P G L (2,\mathbb C)$. Maintenant si $\varphi \in \text{Aut} \Gamma$ il est clair que l'image par $\varphi$ d'une corde de $\Gamma$ est une corde de $\Gamma$ ; ainsi $\varphi \in \text{Aut} \mathcal F$ $\square$

On note que $\text{Aut} \mathcal F$ agit transitivement sur $\mathbb C\mathbb P (3)-\Gamma$ et sur $\Gamma$.
\vskip 3mm
Lorsque $\mathcal F$ est le feuilletage radial en un point 0, il est clair que $\text{Aut} \mathcal F$ est exactement le sous groupe des automorphismes de $\mathbb C\mathbb P(3)$ qui fixent 0.

Examinons le cas o\`u $\mathcal F$ est radial dans les pages d'un livre ouvert $z_{1}/z_{2}=cste$. Commen\c{c}ons par le cas d\'eg\'en\'er\'e o\`u le lieu singulier de $\mathcal F$ se r\'eduit \`a l'axe $z_{1}=z_{2}=0$. On sait qu'alors $\mathcal F$ a deux int\'egrales premi\`eres de base : $z_{2}/z_{1}$ et $\frac{z_{3}-r(\frac{z_{2}}{z_{1}})}{z_{1}}$ o\`u $r$ est une certaine fonction rationnelle non constante. Si l'on \'eclate l'axe $z_{1}=z_{2}=0$ dans $\mathbb C\mathbb P(3)$, on obtient un feuilletage $\widetilde{\mathcal F}$ sur la vari\'et\'e \'eclat\'ee $\mathbb C\widetilde{\mathbb P(3)}$ dont le lieu singulier est donn\'e dans la carte $(z_{1}, t=\frac{z_{2}}{z_{1}}, z_{3})$ par :
$$\text{Sing} \widetilde{\mathcal F}=\{z_{1}=0, z_{3}=r(t) \}$$
Notons que si $\varphi\in \text{Aut} \mathcal F$, alors $\varphi$ pr\'eserve l'axe $z_{1}=z_{2}=0$ ; on peut donc relever $\varphi$ \`a $\mathbb C\widetilde{\mathbb P (3)}$ et obtenir un biholomorphisme $\widetilde{\varphi}\in \text{Aut} \widetilde{\mathcal F}$. Visiblement $\widetilde{\varphi}$ pr\'eserve le diviseur exceptionnel  $(z_{1}=0)$ et laisse invariant  $\text{Sing} \widetilde{\mathcal F}$. De la m\^eme fa\c{c}on si $X$ est un champ de vecteurs sur $\mathbb C\mathbb P(3)$ dont le flot $\varphi_{s}$ est dans $\text{Aut} \mathcal F$ alors $X$ s'\'ecrit :
$$(a_{1} z_{1}+a_{2} z_{2}) \frac{\partial}{\partial z_{1}}+(b_{1} z_{1}+b_{2} z_{2}) \frac{\partial}{\partial z_{2}}+(c_{0} +c_{1}z_{1}+c_{2} z_{2}+c_{3} z_{3}) \frac{\partial}{\partial z_{3}}\leqno (52)$$
$$+(A_{1} z_{1}+A_{2} z_{2}+A_{3} z_{3}) \left(z_{1} \frac{\partial}{\partial z_{1}} +z_{2}\frac{\partial}{\partial z_{2}}+z_{3} \frac{\partial}{\partial z_{3}}\right).$$
La restriction de  l'\'eclat\'e $\widetilde X$ de $X$ au diviseur exceptionnel  $z_{1}=0$ s'\'ecrit~ :
$$(b_{1}+(b_{2}-a_{1}) t-a_{2}t^2) \frac{\partial}{\partial t}+(c_{0}+c_{3}z_{3}+A_{3}z^2_{3})\frac{\partial}{\partial z_{3}} \leqno (53)$$
et doit \^etre tangent \`a $\text{Sing} \widetilde{\mathcal F}$. Ce sera \'evidemment le cas si (53) est identiquement nul. Dans cette \'eventualit\'e $X$ est du type :
$$a_{1}
\left(z_{1}\frac{\partial}{\partial z_{1}}+z_{2}\frac{\partial}{\partial z_{2}}\right)+ (c_{1} z_{1}+c_{2}z_{2}) \frac{\partial}{\partial z_{3}}+ \leqno (54)$$
$$(A_{1} z_{1}+A_{2}z_{2}).\left(z_{1}\frac{\partial}{\partial z_{1}}+ z_{2} \frac{\partial}{\partial z_{2}}+z_{3}\frac{\partial}{\partial z_{3}}\right), $$
et son \'eclat\'e $\widetilde{X}$ s'\'ecrit :
$$\widetilde X=z_{1}\left\{a_{1}+A_{1}z_{1} +A_{2}z_{1}t) \frac{\partial}{\partial z_{1}} +(A_{1}+c_{1}+t(c_{2}+A_{2} z_{3}))\frac{\partial}{\partial z_{3}}\right\}. \leqno (55)$$
Comme $\widetilde X$ doit respecter $\widetilde{\mathcal F}$, dans chaque plan g\'en\'erique on a :
$$[R_{1}, \widetilde X]=h. R_{1}\quad \hbox{\rm o\`u}\quad R_{1}=z_{1} \frac{\partial}{\partial z_{1}}+(z_{3}-r(t)) \frac{\partial}{\partial z_{3}}\leqno (55)$$
d\'efinit $\widetilde{\mathcal F}$ et $h$ est un polyn\^ome en $z_{1}, z_{2}$ \`a param\`etre $t$. En \'ecrivant explicitement (55) on constate que $A_{1}= A_{2}=0$.

Par contre tous les flots des champs :
$$a_{1}\left(z_{1} \frac{\partial}{\partial z_{1}}+z_{2}\frac{\partial}{\partial z_{2}}\right) + (c_{1}z_{1}+c_{2}z_{2}) \frac{\partial}{\partial z_{3}}$$
laissent invariant $\mathcal F$.

Lorsque la fonction rationnelle $r$ est suffisamment g\'en\'erique la composante neutre $\text{Aut}_{0} \mathcal F$ de $\text{Aut} \mathcal F$ se limite au groupe r\'esoluble de dimension 3 : 
$$\text{exp} \left\{a_{1} \left(z_{1}\frac{\partial}{\partial z_{1}}+z_{2} \frac{\partial}{\partial z_{2}}  \right)+(c_{1} z_{1}+ c_{2}z_{2}) \frac{\partial}{\partial z_{3}}\right\} \leqno (56)$$
c'est ce que nous allons voir maintenant tout en classifiant les cas exceptionnels. Consid\'erons un champ $X$ (52) dont l'\'eclat\'e $\widetilde X$ est non identiquement nul sur le diviseur exceptionnel ; alors $\text{Sing} \widetilde{\mathcal F}$ est invariant par le champ donn\'e par (53). Dit autrement la fonction rationnelle $r$ est solution de l'\'equation diff\'erentielle de Riccati :
$$(b_{1}+(b_{2}-a_{1}) t-a_{2}t^2) y'-(c_{0}+c_{3} y+A_{3} y^2)=0.\leqno (57)$$
Ce qui prouve l'affirmation ci-dessus. L'\'equation (57) poss\`ede les solutions constantes $r_{1}$ et $r_{2}$, racines du trin\^ome
$$c_{0}+c_{0}y +A_{3}y^2=0. \leqno (58)$$
Notons que le champ (53) ne peut s'annuler sur $\text{Sing} \mathcal F$ sans \^etre identiquement nul, cas trait\'e pr\'ec\'edemment. Si (58) a deux racines distinctes on peut supposer qu'elles sont en 0 et $\infty$. A changement de coordonn\'ees lin\'eaires pr\`es (57) prend alors l'une des deux formes suivantes :
$$t^2 y'+y=0 \leqno (59)$$
$$t y'+\lambda y=0. \leqno (60)$$
Mais (59) n'a d'autre solutions rationnelles que $y=0$, cas exclus tandis que (60) a des solutions rationnelles $y=r(t)$ non constantes si et seulement si $\lambda \in \mathbb Z$, auquel cas toujours \`a conjugaison pr\`es, de telles solutions sont de type $r(t)=t^n, \ \ n\in \mathbb N$.  Si (58) a une racine double (le cas $A_{3}=0$ se ram\`ene au pr\'ec\'edent) (57) \'equivaut \`a
$$a(t) y'+y^2=0,\ \ \hbox{\rm o\`u\ } a \ \ \hbox{\rm est\ un\ polyn\^ome\  de\ degr\'e\ deux}. \leqno (61)$$
Le seul cas o\`u (61) poss\`ede une solution rationnelle $r(t)$ non constante est celui o\`u $a(t)$ \`a une racine double ; dans cette \'eventualit\'e (61) est conjugu\'ee \`a :
$$t^2 y'- y^2 =0. \leqno (62)$$
qui s'int\`egre en $y(t)=\frac{t}{1+\mu t},\quad \mu \in \mathbb C$.

La fonction rationnelle $r$ est ici l'un des $\frac{t}{1+\mu t}$ et quitte \`a composer par $\frac{t}{1-\mu t}$ on supposera que $r(t)=t$. Nous \'etudions donc le groupe $\text{Aut}\mathcal F_{n}$ o\`u $\mathcal F_{n}$ est le feuilletage en droites associ\'e \`a la fonction $t^n$. En \'ecrivant l'invariance de $\text{Sing} \widetilde{\mathcal F}$ par (53) on obtient :
\vskip 3mm
$\text{Aut}_{0}\mathcal F_{1}=\text{exp} \left\{(a_{1}z_{1}+a_{2}z_{2})\frac{\partial}{\partial z_{1}} + (b_{1} z_{1}+b_{2}z_{2})\frac{\partial}{\partial z_{2}} +(b_{1}+c_{1}z_{1}+c_{2}z_{2}+(b_{2}-a_{1})z_{3})\right.$

$\frac{\partial}{\partial z_{3}}\left.+(A_{1}z_{1}+A_{2}z_{2}-a_{2}z_{3}) \left(z_{1}\frac{\partial}{\partial z_{1}}+z_{2}\frac{\partial}{\partial z_{2}}+ z_{3}\frac{\partial}{\partial z_{3}}\right) a_{i}, b_{i}, c_{i}, A_{i}\in \mathbb C\right\}$
\vskip 3mm
\noindent et pour $n\geq 2$ :
\vskip 2mm
$\text{Aut}_{0}\mathcal F_{n}=\text{exp} \left\{(a_{1}z_{1}+a_{2}z_{2})\frac{\partial}{\partial z_{1}}+\left(b_{2}z_{2}\frac{\partial}{\partial z_{2}}+c_{0}+c_{1}z_{1}+c_{2}z_{2}+n(b_{2}-a_{1}) z_{3}\right) \frac{\partial}{\partial z_{3}}\right.$

$\left.+(A_{1}z_{1}+A_{2}z_{2}+A_{3}z_{3}) \left(z_{1}\frac{\partial}{\partial z_{1}}+z_{2}\frac{\partial}{\partial z_{2}} + z_{3}\frac{\partial}{\partial z_{3}}\right)\right\}.$
\vskip 2mm
Les $\text{Aut}_{0} \mathcal F_{n}$ sont tous de dimension 9 et agissent transitivement sur 
$\mathbb C\mathbb P(3)-\overline{(z_{1}=z_{2}=0}$ 
\vskip 3mm
Venons en  au cas g\'en\'erique d'un feuilletage en droites $\mathcal F$ dans les pages d'un livre ouvert, cas o\`u le lieu singulier $\text{Sing} \mathcal F$ se compose de deux courbes rationnelles  $(z_{1}=z_{2}=0)$ et la courbe  $t\rightsquigarrow (r_{1}(t), t r_{1}(t), r_{3}(t))$ $=R(t)$. Bien s\^ur pour $R$ g\'en\'erique $\text{Aut}_{0}\mathcal F$ est trivial puisque l'image de $R$ en g\'en\'eral n'est pas trajectoire d'un champ de vecteur sur $\mathbb C\mathbb P(3)$. La description des $\text{Aut}_{0}\mathcal F$ dans le cas contraire est "zoologique". L'exemple le plus simple est le suivant.

Supposons que $t\rightsquigarrow R(t)$ soit affine, ie param\`etre une droite g\'en\'erale $D$. On peut supposer que $D$ est donn\'ee par $z_{1}=1, z_{3}=0$. Alors $\text{Aut} \mathcal F$ coincide avec le sous groupe de $\text{Aut} \mathbb C\mathbb P(3)$ qui laisse invariant $\text{Sing} \mathcal F$. Il est isomorphe au projectivis\'e du groupe $G\subset G L (4,\mathbb C)$, o\`u $G$ est engendr\'e par les matrices :

$$\left(
\begin{array}{c|c}
  A  & 0  \\
 \hline  0 &  B 
\end{array}
\right) = A, B\in G L (2,\mathbb C)\quad \hbox{\rm et}\quad 
\left(
\begin{array}{cccc}
 0 &  0 &1&0   \\
  0&0   &  0&1 \\
  1&0   &0&0\\
  0&1&0&0   
\end{array}
\right)
.$$
Il agit encore transitivement sur $\mathbb C\mathbb P(3)-\text{Sing} \mathcal F$.
\section{Remarques et probl\`emes}
La description des solutions rationnelles de l'\'equation aux eikonales en dimension 3 donne bien s\^ur celle des \'equations diff\'erentielles de type~ :
$$Q\left(\frac{\partial f}{\partial z_{1}} , \ \frac{\partial f}{\partial z_{2}},\ \frac{\partial f}{\partial z_{3}} \right)=0,$$
o\`u $Q$ est une forme quadratique de rang maximum. Ces solutions ne sont jamais maximales.

On peut s'interroger sur les solutions Liouvilliennes de l'\'equation au eikonales. 
C'est tout \`a fait naturel puisque les $c. \sqrt{z_{1}^2+z_{2}^2+z_{3}^2}$ en sont des solutions. 
De m\^eme on peut s'int\'eresser aux solutions globales r\'eelles de classe $C^2$ ; on conjecture qu'elles sont affines. La classification des champs de droites en dimension sup\'erieure \`a 4 et ses cons\'equences sur l'\'equation des eikonales restent ouvertes ; soit $P$ un polyn\^ome irr\'eductible tel qu'il existe une solution maximale $f$ rationnelle de $P\left(\frac{\partial f}{\partial z}\right)=0$. On sait que $P=0$ est unirationnelle ; il s'agit de classifier de tels polyn\^omes $P$.
\vskip 3mm
Remerciements \`a Frank Loray et Marie-Annick Paulmier

\end{document}